\documentclass[a4paper,12pt]{amsart}
\usepackage{amscd}       
\usepackage{amsfonts}
\usepackage{amsmath}
\usepackage{amssymb}
\usepackage{amstext}
\usepackage{amsthm}       
\usepackage{xspace}
\usepackage[all]{xy}

\usepackage{graphicx}
\usepackage{url} 

\newcommand{\C}{\mathbb C}
\newcommand{\R}{\mathbb R}
\newcommand{\N}{\mathbb N}
\newcommand{\Q}{\mathbb Q}
\newcommand{\Z}{\mathbb Z} 

\newtheorem{theorem}{Theorem}[section]
\newtheorem{dfn}[theorem]{Definition}
\newtheorem{lemma}[theorem]{Lemma}
\newtheorem{cor}[theorem]{Corollary}
\newtheorem{rem}[theorem]{Remark}
\newtheorem{prop}[theorem]{Proposition}

\begin{document}

\author{Paolo Ghiggini} 
\address{Dipartimento di Matematica \\ Universit\`a di Pisa \\ largo Pontecorvo
5 \\ I--56127 Pisa \\ Italy}
\email{ghiggini@mail.dm.unipi.it}
\title{Ozsv{\'a}th-Szab{\'o} Invariants and fillability of contact structures} 
\date{}
\thanks{The author is a member of EDGE, Research 
Training Network HPRN-CT-2000-00101, supported by The European Human Potential
 Programme.}

\begin{abstract} 
Recently Ozsv{\'a}th and Szab{\'o} defined an invariant of contact 
structures
with values in the Heegaard-Floer homology groups. They also proved that a 
version of the invariant with twisted coefficients is 
non trivial for weakly symplectically fillable contact structures. In this 
article we show that their non vanishing
result does not hold in general for the contact invariant with untwisted 
coefficients. As a 
consequence of this fact Heegaard-Floer theory can distinguish
between  weakly and strongly symplectically fillable contact structures.
\end{abstract} 

\maketitle

\section{Introduction}
Recently Ozsv{\'a}th and Szab{\'o} showed how to associate to any contact
manifold $(Y, \xi)$ an isotopy invariant $c(\xi) \in \widehat{HF}(- Y)/ \pm 1$ in the 
Heegaard-Floer homology of $-Y$ reduced modulo $\pm 1$. They also proved that 
$c(\xi)=0$ if $\xi$ is an overtwisted contact structure, and $c(\xi)$ is a primitive 
element of $\widehat{HF}(-Y) / \pm 1$ if $\xi$ is Stein fillable, \cite{O-Sz:cont}.
One can get rid of the sign indeterminacy in the definition of $c(\xi)$ by
working with the Heegaard--Floer homology with coefficients in $\Z / 2 \Z$.
This is the choice we will do throughout this article. 
The Ozsv{\'a}th-Szab{\'o} contact invariant has already been useful in proving 
tightness of contact structures which resisted to all previously known 
techniques: see for example 
\cite{lisca-stipsicz:3,lisca-stipsicz:4,lisca-stipsicz:5} 

In this article we study the relation between the Ozsv{\'a}th-Szab{\'o} contact 
invariant and the symplectic fillability of contact structures. 
There are two different notions
of symplectic fillability. A contact manifold $(Y, \xi)$ is said to be 
{\em weakly symplectically fillable} if $Y$ oriented by $\xi$ is the 
oriented boundary of a symplectic $4$-manifold $(X, \omega)$ such that $\omega|_{\xi} >0$. 
A contact manifold $(Y, \xi)$ is said to be 
{\em strongly symplectically fillable} if $\xi$ is the 
kernel of a $1$-form $\alpha$ such that $d \alpha = \omega |_Y$. Strong fillability implies 
weak fillability, but the converse is not true. The first example of 
a weakly but not strongly fillable contact manifold was discovered on $T^3$ by
Eliashberg \cite{eliashberg:5}, and more examples were constructed by
Ding and Geiges \cite{ding-geiges:1} on torus bundles over $S^1$ building on 
Eliashberg's. 

We will construct infinitely many weakly fillable contact structures whose 
contact invariant is trivial. These are the first examples of  tight 
contact structures with vanishing Ozsv\'ath--Szab\'o invariant over $\Z / 2 \Z$. 
More precisely, let 
$$M_0 = T^2 \times [0,1] / (\mathbf{v}, 1) = (A \mathbf{v}, 0)$$
be the mapping torus of the map $A:T^2 \to T^2$ induced by the matrix
$\left ( \begin{matrix} 1 & 1 \\
                       -1 & 0
\end{matrix} \right )$. Giroux constructed a family of weakly symplectically
fillable contact structures $\xi_n$ on $M_0$ for $n \in \N^{+}$ as follows. 
Put coordinates $(x,y,t)$ on $T^2 \times \R$ and fix a function $\phi : \R \to \R$. 
For any $n>0$ the $1$-form 
$$\alpha_n = \sin(\phi(t))dx+ \cos(\phi(t))dy$$
on $T^2 \times \R$ defines a contact structure $\xi_n$ on $M_0$  provided that 
\begin{enumerate}
\item $\phi'(t)>0$ for any $t \in \R$
\item $\alpha_n$ is invariant under the action $({\mathbf v},t) \mapsto (A {\mathbf v}, 
t-1)$
\item $(2n-1) \pi \leq \sup \limits_{t \in \R} (\phi(t+1)- \phi(t))< 2n \pi$.
\end{enumerate}
  
The main result of this article is the following theorem.
\begin{theorem}\label{principale}
If $n$ is even, then the Ozsv\'ath--Szab\'o contact invariant $c(\xi_n)$ is trivial.
\end{theorem}

Theorem \ref{principale} should be contrasted with a recent non vanishing 
result for the contact invariant with twisted coefficients proved by Ozsv\'ath 
and Szab\'o.
Associated to any module $A$ over the group ring $\Z[H^1(M, \Z)]$ of $H^1(M, \Z)$
 there is a  
Heegaard--Floer homology group ``with twisted coefficients''  
$\underline{\widehat{HF}}(M; A)$. The ordinary ``untwisted'' Heegaard--Floer 
group is a particular case of this construction with $A= \Z / 2 \Z$. 
See \cite{O-Sz:2}, Section 8.
In this setting the contact invariant $c(\xi)$ can be generalised to an
invariant $\underline{c}(\xi ;A)$ with values in 
$\underline{\widehat{HF}}(-M; A) / \Z[H^1(M, \Z)]^{\times}$, where $\Z[H^1(M, \Z)]^{\times}$ 
denotes the multiplicative group of the invertible elements in $\Z[H^1(M, \Z)]$.

Let $(W, \omega)$ be a weak symplectic filling of the contact manifold $(M, \xi)$.
Following \cite{O-Sz:genus}, we define a $\Z[H^1(M, \Z)]$-module structure on 
$\Z[\R]$ via the ring homomorphism $H^1(M, \Z) \to  \Z[\R]$ defined as
$$\gamma \mapsto T^{\int_M \gamma \land \omega}$$
where $T^r$ denotes the group-ring element associated to the real number $r$.
 The Heegaard-Floer homology group with twisted 
coefficients in the module $\Z[\R]$ will be denoted by 
$\underline{HF}(M; [\omega])$. 
The contact invariant with twisted coefficients of weakly symplectically 
fillable contact structures satisfies the following
non vanishing theorem.

\begin{theorem} \label{bello} (\cite{O-Sz:genus}, Theorem 4.2). 
 Let $(W, \omega)$ be a weak symplectic filling of $(M, \xi)$. Then the associated 
contact invariant $\underline{c}(\xi, [\omega]) \in \underline{\widehat{HF}}(M; [\omega]) /
\Z[H^1(M, \Z)]^{\times}$ is non torsion and primitive.
\end{theorem}

Theorem \ref{bello} implies that the ``untwisted'' Ozsv\'ath--Szab\'o invariant
of a strongly symplectically fillable contact structure is non trivial, 
therefore the contact manifolds $(M_0, \xi_n)$ are not strongly symplectically 
fillable if $n$ is even.
Theorem \ref{principale} shows that, in general, the use of twisted 
coefficients in the non triviality theorem for weakly symplectically
fillable contact structures cannot be avoided, and that the 
Heegaard-Floer theory is subtle enough to 
distinguish between weakly and strongly symplectically fillable contact 
structures. 

\subsection*{Acknowledgements} I thank Ko Honda, Paolo Lisca and Andr{\'a}s 
Stipsicz 
for their encouragement and for many useful discussions. I also thank Peter 
Ozsv{\'a}th for helping me to understand Heegaard-Floer homology and Olga 
Plamenevskaya for answering some questions.
\section{Contact Ozsv\'ath--Szab\'o invariants}
\subsection{Heegaard--Floer homology}
Heegaard--Floer homology is a family of topological quantum field theories
for $Spin^c$ three--manifolds introduced by Ozsv\'ath and Szab\'o in 
\cite{O-Sz:1, O-Sz:2, O-Sz:3}. 
They associate $\Z / 2 \Z$--graded Abelian groups 
$\widehat{HF}(Y, \mathfrak{t})$, $HF^{\infty}(Y, \mathfrak{t})$, 
$HF^-(Y, \mathfrak{t})$, and $HF^+(Y, \mathfrak{t})$ 
to any closed oriented $Spin^c$ $3$--manifold $(Y, \mathfrak{t})$, and  
homomorphisms 
$$F^{\circ}_{W, \mathfrak{s}}: HF^{\circ}(M, \mathfrak{t_1}) \to HF^{\circ}(M, \mathfrak{t_2})$$
to any oriented $Spin^c$ cobordism $(W, \mathfrak{s})$ between two $Spin^c$
manifolds $(M, \mathfrak{t_1})$ and $(M, \mathfrak{t_2})$. Here $HF^{\circ}$ denotes
any of the four functors $\widehat{HF}$, $HF^+$, $HF^-$, and $HF^{\infty}$.
We write $HF^{\circ}(Y)$ for the direct sum $\bigoplus \limits_{\mathfrak{t} \in Spin^c(Y)} 
HF^{\circ}(Y)$ and $F^{\circ}_W$ for the sum $\sum \limits_{\mathfrak{s} \in Spin^c (W)} 
F^{\circ}_{W,\mathfrak{s}}$. $F^{\circ}_W$ is a well defined map because $F^{\circ}_{W,\mathfrak{s}} \neq 0$ only 
for finitely many $Spin^c$--structures on $W$. The homomorphisms between 
Heegaard--Floer homology groups satisfy the following composition rule.
\begin{theorem} (\cite{O-Sz:3}, Theorem 3.4).
Let $(W_1, \mathfrak{s}_1)$ be a $Spin^c$--cobordism between $(Y_1, \mathfrak{t}_1)$
and $(Y_2, \mathfrak{t}_2)$, and let  $(W_2, \mathfrak{s}_2)$ be a 
$Spin^c$--cobordism between $(Y_2, \mathfrak{t}_2)$ and $(Y_3, \mathfrak{t}_3)$.
Denote by $W$ the cobordism between $Y_1$ and $Y_2$ obtained by gluing $W_1$ and 
$W_2$ along $Y_2$. Then 
$$F^{\circ}_{W_2, \mathfrak{s}_2} \circ F^{\circ}_{W_1, \mathfrak{s}_1} = 
\sum \limits_{\begin{matrix} \mathfrak{s} \in Spin^c(W) \\ \mathfrak{s}|_{W_i}=
 \mathfrak{s}_i \end{matrix}}  F^{\circ}_{W, \mathfrak{s}}. $$ 
\end{theorem}
%
The groups $HF^{\circ}(Y, \mathfrak{t})$  are linked to each other by the exact 
triangles

\begin{align}
\xymatrix{ 
& HF^-(Y, \mathfrak{t}) \ar[r] & HF^{\infty}(Y, \mathfrak{t}) \ar[r] & 
HF^+(Y, \mathfrak{t}) \ar `r[d] `[l] `[lll] `[l] [ll] \\   
& & & 
}
\\
\xymatrix{ 
& \widehat{HF}(Y, \mathfrak{t}) \ar[r] & HF^+(Y, \mathfrak{t}) \ar[r] & 
HF^+(Y, \mathfrak{t}) \ar `r[d] `[l] `[lll] `[l] [ll] \\
& & &       
}
\end{align}
These exact triangles are natural, in the sense that they commute with
the maps induced by cobordisms. 

The Heegaard-Floer homology groups $HF^{\circ}(Y, \mathfrak{t})$ have a natural 
$\Z / \mbox{div}(\mathfrak{t})$ 
relative grading, where $\mbox{div}(\mathfrak{t})$ is the divisibility of
$c_1(\mathfrak{t})$ in $H^2(Y, \Z)$. it was shown in \cite{O-Sz:4} that, when 
$c_1(\mathfrak{t})$ is a torsion
element, the
relative $\Z$--grading admits a natural lift to an absolute
$\Q$--grading.  In conclusion, for a torsion $Spin^c$--structure $\mathfrak{t}$ 
on $Y$ the Ozsv\'ath--Szab\'o homology groups $HF^{\circ} (Y, \mathfrak{t})$ split as
$$HF^{\circ} (Y, \mathfrak{t})=\bigoplus_{d \in  \Q} HF^{\circ}_d (Y, \mathfrak{t}).$$ 
When $\mathfrak{t} \in Spin^c(Y)$ has torsion first Chern class, there is an 
isomorphism between the homology groups $\widehat{HF}_d (Y,\mathfrak{t})$ and 
$\widehat{HF}_{-d}(-Y,\mathfrak{t})$.

\begin{prop} (See \cite{O-Sz:4}, Theorem 7.1).
Let $(W, \mathfrak{s} )$ be a $Spin^c$ cobordism between two $Spin^c$ manifolds
$(Y_1,\mathfrak{t}_1) $ and $(Y_2, \mathfrak{t}_2)$. If the $Spin^c$ structures 
$\mathfrak{t}_i$ have both torsion first Chern class and $x \in HF^{\circ} 
(Y_1, \mathfrak{t}_1)$ is a homogeneous element of degree $d(x)$, 
then $F_{W, \mathfrak{s}}(x) \in HF^{\circ} (Y_2, \mathfrak{t}_2)$ is also
homogeneous of degree
$$d(x) + \frac{1}{4}(c_1^2(\mathfrak{s})-3 \sigma (W)-2 \chi (W)).$$
\end{prop}
Notice that $F^{\circ}_W$  might map a homogeneous element $x\in HF^{\circ}_d (Y_1, 
\mathfrak{t}_1)$ into a non homogeneous element $F^{\circ}_W(x) \in HF^{\circ}(Y_2)$.

\subsection{Definition of the contact invariants.}
The Ozsv\'ath--Szab\'o contact invariant is defined using the correspondence 
between contact structures and open book decompositions of three--manifolds
recently discovered by Giroux. An {\em open book decomposition} of a
$3$--manifold $Y$ is a fibred link $B \subset Y$ together with a fibration
$\pi : Y \setminus B \to S^1$. The link $B$ is called the {\em binding} of the open book 
decomposition and the union of a fibre of $\pi : Y \setminus B \to S^1$ with $B$ is called
a {\em page}.
\begin{dfn} (\cite{giroux-obd}, Definition 1). 
Let $(Y, \xi)$ be a contact $3$--manifold. An open book decomposition $(B, \pi)$
of $Y$ is said to be 
adapted to $\xi$ if:
\begin{enumerate}
\item $B$ is transverse to $\xi$,
\item $\xi$ is defined by a contact form $\alpha$ such that $d \alpha$ is a symplectic 
form on any fibre of $\pi$,
\item the orientation of $B$ induced by the contact structure coincides with
the orientation as boundary of the fibres of $\pi$ oriented by $d \alpha$.
\end{enumerate}
\end{dfn}

By \cite{giroux-obd} Theorem 3 any contact structure on a 
three manifold admits an adapted open book decomposition. This open book 
decomposition is not unique, in fact two open book decompositions which differ
by the positive plumbing  of an annulus are adapted to isotopic contact
structures. See \cite{giroux-obd} Section B. After positive plumbing, we
can assume that the binding is connected and pages have genus $g \geq 2$. Adding 
a $2$--handle along $B$ with the 
framing induced by a page we form a cobordism $V$ between $Y$ and $Y_0$, where 
$Y_0$ is a $3$-manifold fibred over $S^1$ with fibres of genus $g \geq 2$. On $Y_0$
there is a canonical $Spin^c$--structure $\mathfrak{t}_0$ induced by the 
fibration. $\widehat{HF}(-Y_0, \mathfrak{t}_0)= \Z /2 \Z \oplus \Z /2 \Z$ with the 
summands lying in different degrees for the absolute $\Z /2 \Z$ grading, while 
$HF^+(-Y_0, \mathfrak{t}_0)= \Z /2 \Z$. See 
\cite{O-Sz:cont} Section 3. We fix a distinguished element $c_0 \in 
\widehat{HF}(-Y_0, \mathfrak{t}_0)$ as the homogeneous element of 
$\widehat{HF}(-Y_0, \mathfrak{t}_0)$ which is mapped to the non zero element of
$HF^+(-Y_0, \mathfrak{t}_0)$ by the natural map $\widehat{HF}(-Y_0, \mathfrak{t}_0)
\to HF^+(-Y_0, \mathfrak{t}_0)$. We denote by $\overline{V}$ the cobordism $V$ 
turned upside--down, so that $\overline{V}$ is a cobordism between 
$-Y_0$ and $-Y$. 
\begin{dfn}
The {\em Ozsv\'ath--Szab\'o contact invariant} of a contact $3$--manifold
$(Y, \xi)$ is the element $c(\xi) \in \widehat{HF}(-Y)$ defined by
$$c(\xi)= \widehat{F}_{\overline{V}}(c_0).$$
\end{dfn} 
 By \cite{O-Sz:cont} Theorem 1.3 $c(\xi)$ is independent of the choice of the
open book decomposition adapted to $\xi$ and is an isotopy invariant.
The Ozsv\'ath--Szab\'o contact invariant is non trivial and 
detects important topological properties of the contact structures, in fact
\begin{theorem} (\cite{O-Sz:cont}, Theorem 1.4 and Theorem 1.5).
If $(Y, \xi)$ is overtwisted, than $c(\xi)=0$. If $(Y, \xi)$ is Stein fillable, then
$c(\xi) \neq 0$.
\end{theorem}
The Ozsv\'ath--Szab\'o contact invariant $c(\xi)$ encodes the homotopy invariants 
of $\xi$, see \cite{O-Sz:cont}, Proposition 4.6.
Any contact structure $\xi$ on a $3$--manifold 
$Y$ determines a $Spin^c$--structure $\mathfrak{t}_{\xi}$ on $Y$, then 
$c(\xi) \in \widehat{HF}(-Y, \mathfrak{t}_{\xi})$. If the first Chern class of
$\xi$ is torsion, by \cite{gompf:1} Theorem 4.16 the homotopy type of $\xi$ is 
determined by the $Spin^c$--structure $\mathfrak{t}_{\xi}$ and by the $\Q$--valued 
Gompf invariant $d_3(\xi)$ defined as follows.

\begin{dfn} (See \cite{gompf:1}, Definition 4.2).
Let $\xi$ be an oriented tangent plane field on the $3$--manifold $Y$ with 
torsion first Chern class, and let $(X, J)$ be a almost complex $4$--manifold
such that $Y$ is the boundary of $X$ and $\xi= TY \cap J(TY)$ is the field of 
complex lines in $TY$. Then we define
$$d_3= \frac 14 (c_1(J)^2 -2 \chi(X) -3 \sigma(X))$$
where $\chi$ denote the Euler characteristic, $\sigma$ the signature, and $c_1(J)^2$ is
defined because $c_1(\xi)=c_1(J)|_Y$ is torsion.
\end{dfn} 

By \cite{O-Sz:cont}, Proposition 4.6, if $c_1(\xi)$ is 
a torsion element of $H^2(Y, \Z)$, then $c(\xi)$ is an homogeneous element of 
degree $-d_3(\xi)- \frac 12$. 

\begin{theorem} \label{chirurgia}(\cite{O-Sz:cont}, Theorem 4.2 and 
\cite{lisca-stipsicz:3}, Theorem 2.3).
If the contact manifold $(Y', \xi')$ is obtained from the contact manifold
$(Y, \xi)$ by Legendrian surgery along a 
Legendrian knot $L$, and $W$ is the cobordism between $Y$ and $Y'$ obtained
by adding a $2$--handle to $Y \times [0,1]$ along $L \times \{ 1 \}$ with framing $-1$ 
with respect to the contact framing, then
$$\widehat{F}_{\overline{W}}(c(\xi'))= c(\xi)$$ 
where $\overline{W}$ denotes the cobordism $W$ turned upside--down.
\end{theorem}
The space of oriented contact structures on $Y$ has a natural involution. 
\begin{dfn}
For any contact structure $\xi$ on a $3$--manifold $Y$ we denote by 
$\overline{\xi}$ the contact structure on $Y$ obtained from $\xi$ by inverting the 
orientation of the planes. 
\end{dfn}

This operation is compatible with the conjugation of the $Spin^c$-structure 
defined by the contact structure, in fact
$\mathfrak{t}_{\overline{\xi}}= \overline{\mathfrak{t}_{\xi}}$. 
There is an isomorphism $\mathfrak{J}: HF^{\circ}(-Y, \mathfrak{s}) \to 
HF^{\circ}(-Y, \overline{\mathfrak{s}})$ defined in \cite{O-Sz:2}, Theorem 2.4. We 
recall that the isomorphism $\mathfrak{J}$ preserves the $\Z /2 \Z$--grading
of the Heegaard--Floer homology groups and is a natural 
transformation in the following sense.
 
\begin{prop}(\cite{O-Sz:3}, Theorem 3.6)
Let $(W, \mathfrak{s})$ be a $Spin^c$-cobordism between 
$(Y_1, \mathfrak{t}_1)$ and $(Y_2, \mathfrak{t}_2)$. Then  
the following  diagram 

$$\begin{CD}
HF^{\circ}(Y_1, \mathfrak{t}_1) & @> F^{\circ}_{W, \mathfrak{s}} >> & HF^{\circ}(Y_2, \mathfrak{t}_2) \\
@VV \mathfrak{J} V &   & @VV \mathfrak{J} V \\
HF^{\circ}(Y_1, \overline{\mathfrak{t}}_1) & @> F^{\circ}_{W, \overline{\mathfrak{s}}} >> &
 HF^{\circ}(Y_2, \overline{\mathfrak{t}}_2)
\end{CD}$$
commutes. 
\end{prop}
The isomorphism $\mathfrak{J}$ commutes also with the maps in the 
exact triangles (1) and (2) relating the different Heegaard--Floer homology
groups.

\begin{theorem}\label{coniugazione}
Let $(Y, \xi)$ be a contact manifold, then 
$$c(\overline{\xi})= \mathfrak{J}(c(\xi)).$$
\end{theorem}
\begin{proof}
If $(B, \pi)$ is an open book decomposition adapted to $\xi$, then the open book 
decomposition $(-B, \overline{\pi})$, where $-B$
denotes the binding $B$ with opposite orientation and $\overline{\pi}$ is the
composition of $\pi$ with the complex conjugation on $S^1$,  is adapted to 
$\overline{\xi}$. The pages of
$(-B, \overline{\pi})$ are the pages of $(B, \pi)$ with opposite orientation, so
the fibration on $Y_0$ induced by $(-B, \overline{\pi})$ differs from the 
fibration induced by
$(B, \pi)$ for the orientation of the fibres, therefore its canonical 
$Spin^c$--structure is the conjugate of $\mathfrak{t}_0$. The commutative
diagram
$$\begin{CD}
\widehat{HF}(-Y_0', \mathfrak{t}_0) & @>>> & HF^+(-Y_0, \mathfrak{t}_0) \\
@V \mathfrak{J} VV & & @V \mathfrak{J} VV \\
\widehat{HF}(-Y_0', \overline{\mathfrak{t}}_0) & @>>> & HF^+(-Y_0, 
\overline{\mathfrak{t}}_0) \\
\end{CD}$$
together with the fact that $\mathfrak{J}$ is an isomorphism and preserves
the $\Z / 2 \Z$--grading of the Heegaard--Floer homology groups shows that
 the distinguished element of $\widehat{HF}(-Y_0', \overline{\mathfrak{t}}_0)$
is $\overline{c}_0= \mathfrak{J}(c_0)$, therefore
$$c(\overline{\xi})= \widehat{F}_{\overline{V}}(\overline{c}_0) = 
\widehat{F}_{\overline{V}}(\mathfrak{J}(c_0)) = \mathfrak{J}(\widehat{F}_{\overline{V}}(c_0))= 
\mathfrak{J} (c(\xi)).$$
\end{proof}

\subsection{Ozsv\'ath--Szab\'o contact invariants of strongly symplectically 
fillable contact structures}

In this section we prove a non vanishing theorem for the Ozsv\'ath--Szab\'o contact
invariant of strongly symplectically fillable contact structures. This theorem 
can be easily derived as a 
corollary of the more general non vanishing Theorem \ref{bello} proved by 
Ozsv\'ath and Szab\'o using the twisted coefficients, however it is also possible 
to  adapt the 
proof of Theorem \ref{bello}, so that we do not need to use Heegaard-Floer
homologies with twisted coefficients. We choose this second option, but the
proof requires some more Heegaard--Floer machinery.

From the exact triangle (1) we define a 
fifth group $HF^{red}(Y, \mathfrak{t})$ as the kernel of the map
$$HF^-(Y, \mathfrak{t}) \to HF^{\infty}(Y, \mathfrak{t})$$
 or, equivalently, as the cokernel of the map 
$$HF^{\infty}(Y, \mathfrak{t}) \to HF^+(Y, \mathfrak{t}).$$
The group $HF^{red}(Y, \mathfrak{t})$ is always finitely generated.
Let $W$ be an oriented cobordism between the $3$--manifolds $Y_1$ and $Y_2$.
An {\em admissible cut} of $W$ 
(\cite{O-Sz:3}, Definition 8.3) is a $3$--manifold
$N \subset W$ which divides $W$ into two pieces $W_1$ and $W_2$ such that  $b_2^+(W_i)>0$
for $i=1,2$, and the connecting homomorphism $\delta: H^1(N, \Z) \to H^2(W, \partial W)$
of the Meyer--Vietoris sequence of the pair $(W_1, W_2)$ is trivial. It is shown 
in \cite{O-Sz:3},
Example 8.4 that an admissible cut of $W$ always exists if $b_2^+(W)>1$. 
By \cite{O-Sz:3} Lemma 8.2 the maps 
\begin{align*}
& F^{\infty}_{W_1, \mathfrak{s}}: HF^{\infty}(Y_1, \mathfrak{s}|_{Y_1}) \to HF^{\infty}(N, 
\mathfrak{s}|_N) \\
& F^{\infty}_{W_1, \mathfrak{s}}: HF^{\infty}(N, \mathfrak{s}|_N) \to HF^{\infty}(Y_2, 
\mathfrak{s}|_{Y_2})
\end{align*}
 vanish for any 
$Spin^c$--structure $\mathfrak{s}$ on $W$, therefore an easy diagram chase on 
the exact triangle (1) allows us to define a ``mixed'' homomorphism 
$F^{mix}_{W, \mathfrak{s}}: HF^-(Y_1, \mathfrak{t}_1) \to  HF^+(Y_2, \mathfrak{t}_2)$
which factors through $HF^{red}(N, \mathfrak{s})$.  By \cite{O-Sz:3}, Theorem 
8.5 the mixed map $F^{mix}_{W, \mathfrak{s}}$ does not depend on the particular 
admissible cut used to define it. 

The mixed map can be used to define a numerical invariant of smooth
four--manifolds with $b_2^+>1$ which is conjecturally equal to the 
Seiberg--Witten invariant. If $X$ is a closed oriented $4$--manifold, after 
removing two 
balls we can view it as a cobordism from $S^3$ to $S^3$. The groups 
$HF^+(S^3)$ and $HF^-(S^3)$ have distinguished elements $\Theta^+$ and $\Theta^-$ which are
the non trivial elements in minimal (resp. maximal) 
degree. See \cite{O-Sz:2}, Section 3 for the computation of the Heegaard--Floer
homology groups of $S^3$. The four--dimensional invariant of $X$ is the map
$$\Phi_X: Spin^c(X) \to \Z /2 \Z$$
where $\Phi_X(\mathfrak{s})$ is defined as the coefficient of 
$\Theta^+$ in $F^{mix}_{X, \mathfrak{s}}(\Theta^-)$.  

We denote by $c^+(\xi)$ the image of $c(\xi)$ in $HF^+(-Y)$. 
Theorem \ref{chirurgia} can be refined in the following way. 
\begin{lemma} \label{canonica} Suppose that $(Y', \xi')$ is obtained from 
$(Y, \xi)$ by Legendrian surgery on a Legendrian link $L$ 
and that $(W, \omega)$ is the symplectic cobordism from $(Y, \xi)$ to $(Y', \xi')$
induced by this surgery. Then we have
$$F^+_{\overline{W}, \mathfrak{k}}(c^+(\xi'))=c^+(\xi)$$
for the canonical $\mbox{spin}^c$-structure $\mathfrak{k}$ associated to the 
symplectic structure on $W$, and  
$$F^+_{\overline{W}, \mathfrak{s}}(c^+(\xi'))=0$$
 for any $\mbox{spin}^c$-structure $\mathfrak{s}$ on $W$ with $\mathfrak{s} \neq 
\mathfrak{k}$.
\end{lemma}
\begin{proof}
As in the proof of \cite{lisca-stipsicz:3} Theorem 2.3 there exists an open 
book decomposition of $Y$ adapted to the contact structure $\xi$ so that 
the surgery link lies on a page. We can also assume that the binding is 
connected and the pages have genus $g>1$.
An open book decomposition adapted to $\xi'$ is obtained from the open book 
decomposition adapted to $\xi$ by composing the monodromy with right--handed 
Dehn twists along the surgery link.
Let $Y_0$ and 
$Y_0'$ be the $3$-manifolds obtained from $Y$ 
and $Y'$ respectively by $0$-surgery on the binding, and let $V$, $V'$ be 
the induced cobordisms. The surgery on $L$ induces cobordisms $W$ between $Y$
and $Y'$ and  $W_0$ from $Y_0$ to $Y_0'$. Both $Y_0$ and $Y_0'$ are surface bundles 
over $S^1$, and $W_0$ admits a Lefschetz fibration over the annulus. Let 
$\mathfrak{t}_0$ and $\mathfrak{t}_0'$ be the $\mbox{Spin}^c$-structures on $Y_0$ 
and $Y_0'$ respectively determined by the fibration, and let $\mathfrak{k}_0$ 
be the canonical $\mbox{Spin}^c$-structure on $W_0$ determined
by the Lefschetz fibration. By \cite{O-Sz:symp}, Theorem 5.3, 
$$F^+_{\overline{W}_0, \mathfrak{k}_0}: HF^+(-Y_0', \mathfrak{t}_0') \to HF^+(-Y_0, 
\mathfrak{t}_0)$$
is an isomorphism, while the maps
$$F^+_{\overline{W}_0, \mathfrak{s}}: HF^+(-Y_0', \mathfrak{t}') \to HF^+(-Y_0, 
\mathfrak{t})$$
are trivial when $\mathfrak{s} \neq \mathfrak{k}_0$.

Let $W'$ be the cobordism $W'= W_0 \cup_{Y_0'} V' = V \cup_{Y}W$ from $Y_0$ to $Y'$. 
Since the cobordism $V'$ is obtained by adding a unique $2$-handle along a 
homologically non trivial curve, the restriction map $H^2(W', \Z) \to 
H^2(W_0, \Z)$ is an isomorphism, therefore there is a unique 
$\mbox{Spin}^c$-structure $\mathfrak{k}_0'$ on $W$ which extends 
$\mathfrak{k}_0$. By the composition formula \cite{O-Sz:3} Theorem
3.4 $F^+_{W', \mathfrak{k}_0'}= F^+_{V'} \circ F^+_{W_0, \mathfrak{k}_0}$ 
and for any other  $\mbox{Spin}^c$-structure 
$\mathfrak{s} \neq \mathfrak{k}_0'$ the map $F^+_{X, \mathfrak{s}}$ is trivial.
Let $\mathfrak{s}'$ be the restriction of $\mathfrak{k}_0'$ to $W$, then the 
diagram
$$\begin{CD}
HF^+(-Y_0', \mathfrak{t}_0') & @> F^+_{\overline{W}_0, \mathfrak{k}_0} >> & HF^+(-Y_0, 
\mathfrak{t}_0) \\
@V F^+_{\overline{V}'} VV & & @V F^+_{\overline{V}} VV \\ 
HF^+(-Y', \mathfrak{t}_{\xi'}) & @> F^+_{\overline{W}, \mathfrak{s}'} >> & HF^+(-Y, \mathfrak{t}_{\xi}) 
\end{CD}$$
commutes and $F^+_{W, \mathfrak{s}} = 0$ for any $\mathfrak{s} \neq 
\mathfrak{s}'$.
To finish the proof, we have to identify $\mathfrak{s}'$ with $\mathfrak{k}$.

By \cite{eliashberg-fill}, Theorem 1.1, the symplectic structure induced by 
the Lefschetz fibration on $W_0$ extends over the $2$-handle $V'$, thus we 
obtain a symplectic structure $\omega'$ on $W'$ with canonical 
$\mbox{Spin}^c$-structure $\mathfrak{k}_0'$. The restriction of $\omega'$ to $W$ 
coincides with the symplectic structure on $W$ induced by the Legendrian 
surgery, therefore $\mathfrak{s}' = \mathfrak{k}$.
\end{proof}

We have stated Lemma \ref{canonica} in the 
form in which we are going to use it, however it can be proved in the same way
for the stronger contact invariant in $\widehat{HF}(-Y)$ with integer 
coefficients. 

\begin{lemma} \label{piccolo}
Let $(Y, \xi)$ be a contact manifold, then 
there exists a concave symplectic filling $(W', \omega_{W'})$ of $(Y, \xi)$ with 
canonical $Spin^c$-structure $\mathfrak{k}_{W'}$ such that $b_2^+(W')>1$ and 
$$c^+(\xi)=  F^{mix}_{\overline{W'}, \mathfrak{k}_{W'}}(\Theta^-).$$
\end{lemma}

\begin{proof}
Combining \cite{etnyre-honda:4}, Theorem 1.1 and \cite{etnyre-fill} Lemma 3.1  
 there is a Stein fillable contact manifold $(Y', \xi')$ and
a symplectic cobordism $(V_1, \omega_{V_1})$ from $(Y, \xi)$ to $(Y', \xi')$ so that $Y'$ 
is a rational homology sphere and $V_1$ is composed by $2$-handles attached in 
a Legendrian way.  By \cite{plam:1} Lemma 1 there is a 
concave filling $(V_2, \omega_{V_2})$ of $(Y', \xi')$ with canonical 
$Spin^c$-structure $\mathfrak{k}_{V_2}$ such that $b_2^+(V_2)>1$ and
$c^+(\xi')=  F^{mix}_{\overline{V_2}, \mathfrak{k}_{V_2}}(\Theta^-)$. 

Let $(W', \omega_{W'})$ be the concave 
filling of $(Y, \xi)$ obtained by gluing $(V_1, \omega_{V_1})$ and $(V_2, \omega_{V_2})$ along 
$(Y', \xi')$, and let $\mathfrak{k}_{V_1}$ be the 
canonical $Spin^c$-structures of $(V_1, \omega_{V_1})$.
 Since $Y$ is a rational homology sphere, $H^2(W', \Z)=H^2(V_1, \Z) \oplus H^2(V_2, \Z)$
therefore there exists a  unique  $Spin^c$-structure $\mathfrak{k}_{W'}$ on $W'$ 
which restricts to $\mathfrak{k}_{V_1}$ on $V_1$ and to $\mathfrak{k}_{V_2}$ on 
$V_2$. The composition formula \cite{O-Sz:3},
Theorem 3.4, together with Lemma \ref{canonica}, yields 
$$c^+(\xi)=F^+_{\overline{V_1}, \mathfrak{k}_{V_1}} \circ F^{mix}_{\overline{V_2}, 
\mathfrak{k}_{V_2}} (\Theta^-)= F^{mix}_{\overline{W}', \mathfrak{k}_{W'}}(\Theta^-).$$ 
\end{proof}

\begin{theorem} \label{bruciato}
Let $(Y, \xi)$ be a strongly symplectically fillable contact manifold, then
$c(\xi) \neq 0$.
\end{theorem}

\begin{proof}
Let $(W_1, \omega_1)$ be a strong symplectic filling of $(Y, \xi)$, and let $(W_2, \omega_2)$
be the concave symplectic filling considered in Lemma \ref{piccolo}. Gluing
$(W_1, \omega_1)$ and $(W_2, \omega_2)$ we obtain a closed symplectic manifold $(X, \omega)$ with
$b_2^+(X)>1$.
The composition formula \cite{O-Sz:3} Theorem 3.4 gives
$$F^+_{\overline{W}_2, \mathfrak{k}_{W_2}}(c(\xi))= 
F^+_{\overline{W}_2, \mathfrak{k}_{W_2}} \circ F^{mix}_{\overline{W}_1, \mathfrak{k}_{W_1}}
(\Theta^-)=
\sum \limits_{\begin{matrix} \mathfrak{s} \in Spin^c(X) \\ \mathfrak{s}|_{W_i}=
 \mathfrak{k}_{W_i} \end{matrix}}  F^{mix}_{X, \mathfrak{s}} (\Theta^-) = \sum 
\limits_{\begin{matrix} \mathfrak{s} \in Spin^c(X) \\ \mathfrak{s}|_{W_i}=
 \mathfrak{k}_{W_i} \end{matrix}}  \Phi_X(\mathfrak{s}).$$
One of the $Spin^c$--structures in the sum is the canonical one $\mathfrak{k}_X$
coming from the symplectic structure on $X$. For any other 
$Spin^c$--structure $\mathfrak{s}$ in the sum we have
$c_1(\mathfrak{s})-c_1(\mathfrak{k}_X) \in \delta (\alpha(\mathfrak{s}))$ for $\alpha(\mathfrak{s})
\in H^1(Y, \Z)$, where $\delta$ is the 
homomorphism $H^1(Y) \to H^2(X)$ in the Meyer--Vietoris exact sequence for the
pair $(W_1, W_2)$, therefore 
$$\langle c_1(\mathfrak{s})-c_1(\mathfrak{k}_X), [\omega] \rangle_X = \langle \alpha(\mathfrak{s}), [\omega|_Y] \rangle_Y =0$$
in fact $\omega|_Y$ is exact because $W_1$ is a strong filling.

By \cite{O-Sz:symp} Theorem 1.1 the only non zero term in the sum is 
$\Phi_X(\mathfrak{k}_X)=1$, therefore $F^+_{\overline{W}_1, \mathfrak{k}_{W_1}}
(c^+(\xi))= \Theta^+$ which implies that $c^+(\xi) \neq 0$ . In turn, 
this implies that $c(\xi) \neq 0$.
\end{proof}
\begin{rem}
Actually the proof of Theorem \ref{bruciato} proves the stronger fact that, if
we see $(W_1, \omega_1)$ as a symplectic cobordism between the standard $(S^3, \xi_0)$
and $(Y, \xi)$, then 
$$F^+_{\overline{W}_1, \mathfrak{k}_{W_1}}(c(\xi))=c(\xi_0)$$
for the canonical $Spin^c$--structure of $(W_1, \omega_1)$ and
$$F^+_{\overline{W}_1, \mathfrak{s}}(c(\xi))=0$$
for any other $Spin^c$--structure on $(W_1, \omega_1)$ with
$$\langle c_1(\mathfrak{s}), [\omega_1] \rangle_{W_1} = \langle c_1(\mathfrak{k}_{W_1}), [\omega_1] \rangle_{W_1}.$$
\end{rem}

\section{Weakly fillable contact structures with trivial untwisted $\Z /2 \Z$
Ozsv{\'a}th-Szab{\'o} contact invariant}

\subsection{Tight contact structures on $M_0$}
Let $M_0$ be the $T^2$-bundle over  $S^1$ with monodromy map
$A: T^2 \times \{ 1 \} \to T^2 \times \{ 0 \}$ given by 
$A= \left ( \begin{matrix}
                          1 & 1 \\
                         -1 & 0 
                                 \end{matrix} \right )$.
%

Put coordinates $(x,y,t)$ on $T^2 \times \R$ and fix a function $\phi : \R \to \R$. For 
any $n>0$ the $1$-form 
$$\alpha_n = \sin(\phi(t))dx+ \cos(\phi(t))dy$$
on $T^2 \times \R$ defines a contact structure $\xi_n$ on $M_0$  provided 
that 
\begin{enumerate}
\item $\phi'(t)>0$ for any $t \in \R$
\item $\alpha_n$ is invariant under the action $({\mathbf v},t) \mapsto (A {\mathbf v}, 
t-1)$
\item $(2n-1) \pi \leq \sup \limits_{t \in \R} (\phi(t+1)- \phi(t))< 2n \pi$
\end{enumerate}

The main results about this family of contact structures are the following.
\begin{theorem}(\cite{giroux:3}, Proposition 2 and Theorem 6).
The contact structures $\xi_n$ do not depend
 on the function $\phi$ up to isotopy, and are all universally tight and distinct.
\end{theorem}

\begin{theorem}(\cite{honda:2}, Theorem 0.1). \label{honda}
 The tight contact structures $\xi_n$ are the only tight contact structures
on $M_0$ up to isotopy
\end{theorem}

\begin{theorem} (\cite{ding-geiges:1}, Theorem 1). \label{geiges}
For any $n \in \N$, $\xi_n$ is weakly symplectically fillable. There is a  
number $n_0$ such that, for
any $n>n_0$, $\xi_n$  is not strongly symplectically fillable.
\end{theorem}

The fibration on $M_0$ admits a transverse $1$--dimensional foliation induced 
by the foliation by segments on $T^2 \times [0,1]$. Let $F$ be the image of 
$\{ 0 \} \times [0,1]$ in $M_0$, then $F$ is Legendrian 
with respect to the contact structure $\xi_n$ for all $n$. 

The manifold $M_0$ has  a presentation as $0$-surgery on the right-handed
trefoil knot $K$, in fact the complement of $K$ in $S^3$ fibres over 
$S^1$ with fibre the holed torus and the monodromy acts on the homology of the 
fibre as
$A= \left ( \begin{matrix}
                          1 & 1 \\
                         -1 & 0 
                                 \end{matrix} \right )$ for some choice of 
coordinates in the fibre. Moreover the identification between $M_0$ 
and the $0$--surgery on $K$ can be chosen so that the complement of a tubular
neighbourhood of $K$ in $S^3$  is mapped diffeomorphically into 
the complement of a tubular neighbourhood of $F$ in $M_0$ and the 
meridian of $K$ is mapped to a longitude of $F$.

We perform a change of coordinates in a neighbourhood of $F$ to determine
what longitude of $F$ corresponds to the meridian of $K$ and to compute the 
twisting number of $\xi_n$ along $F$ induced by this longitude. 
\begin{lemma} \label{algebra}
Let $R= \left ( \begin{matrix}
                               \frac 12 & \frac{\sqrt{3}}{2} \\
                   - \frac{\sqrt{3}}{2} & \frac 12
\end{matrix} \right )$ be the rotation by angle $- \frac{\pi}{3}$. Then $A$ is
conjugate to $R$ in $GL^+(2,\R)$.
\end{lemma}
\begin{proof}
$A$ and $R$ are conjugated in $GL(2, \C)$ because they have the same 
characteristic polynomial with distinct roots, therefore they are
conjugate also in $GL(2, \R)$ because they are both real. Let $B \in GL(2, \R)$ 
be a matrix such that $BAB^{-1}=R$. For any $x \in \R^2 \setminus \{ 0 \}$ we 
have $x \land Ax \neq 0$ because $A$ has no real eigenvalues, therefore,
after identifying $\bigwedge^2\R^2$ to $\R$ using the canonical basis,
$x \land Ax$ has constant sign as a function $\R^2 \to \R$. A direct computation at
$x= \binom{0}{1}$ shows that $x \land Ax$ is 
negative. For the same reason, $x \land Rx$ is also negative, therefore 
$\det B>0$ because $x \land Rx= B^{-1}Bx \land B^{-1}ABx= (\det B)^{-1} Bx \land ABx$.
\end{proof}

\begin{figure}
\includegraphics[width=6cm]{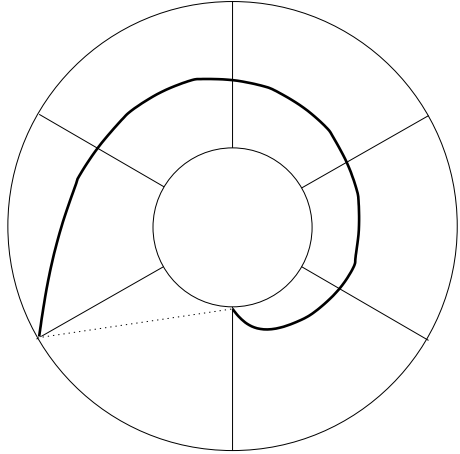}
\caption{The boundary of $V$. The inner circle is glued to the outer one after
a rotation of $- \frac{\pi}{3}$. The dotted line closes to a longitude of $V$, 
the radial lines close to a leaf of the transverse foliation and the bold 
line closes to a 
dividing curve for $\xi_0$.}
\label{figura.fig}
\end{figure}

\begin{lemma}\label{tn}
The twisting number of $\xi_n$ along the Legendrian curve $F$ is
$tn(F, \xi_n)= -n$
\end{lemma}
\begin{proof}
Let $U$ be a small $A$-invariant neighbourhood of $(0,0)$ in $T^2 = \R^2 / \Z^2$
 so that 
$$V= U \times [0,1]/(\mathbf{v}, 1)=(A \mathbf{v}, 0)$$
 is a standard neighbourhood of $F$. Then $B^{-1}$ is defined on $U$ and 
$U_0=B^{-1}(U)$ is a 
$R$-invariant neighbourhood of $(0,0)$, i.~e. a disc centred in $(0,0)$.
In the coordinates $(x',y',t)$ of $U_0 \times \R$ the $1$-form $\alpha_n$ can be written 
 as 
$$\alpha_n= \sin (2 \pi (n+ \frac 56)t)dx'+ \cos (2 \pi (n+ \frac 56)t)dy'.$$
By Lemma \ref{algebra} the leaves of the transverse foliation in the boundary of
the neighbourhood of $K$
have slope $- \frac 16$, therefore they intersect the meridian of $K$ once. 
If we put coordinates $(\theta ,t)$ on $\partial U_0 \times I$, then the longitude of $F$ 
corresponding to the meridian of $K$ is the 
image in $\partial V$ of the arc $t \mapsto (e^{i t \frac{\pi}{3}}, t)$ (the dotted curve in Figure
\ref{figura.fig}) because it intersect the 
leaves of the transverse foliation only once. A dividing curve of 
$\xi_n$ is isotopic to the image of the arc $t \mapsto (e^{-2 \pi (n+ \frac 56)t}, t)$ therefore
 the twisting number of $\xi_n$ along $F$, which is the algebraic intersection 
of a dividing curve with the longitude, is $-n$.  Figure \ref{figura.fig} 
shows what happens for $n=1$.
\end{proof}

\begin{lemma}\label{max}
If $L \subset M_0$ is a Legendrian curve which is smoothly isotopic to $F$, then
$tn(L, \xi_n) \leq tn(F, \xi_n)$
\end{lemma}
\begin{proof}
Since $A^6=I$, $M_0$ has a six-fold cover with total space $T^3$ induced by a cover
of $S^1$. Let $\widehat{F}_3$ and $\widehat{L} \subset T^3$ be the pre-images of 
$F$ and $L$ respectively. By \cite{kanda:1}, Theorem 7.6, $\widehat{F}_3$
maximises the twisting number in its smooth isotopy class. The lemma follows 
from the obvious monotonicity of the twisting number under finite coverings.
\end{proof}

Since the right-handed trefoil can be put in Legendrian form
with Thurston-Bennequin invariant $1$, this surgery presentation yields a 
Stein fillable contact structure on $M_0$.
\begin{prop}\label{identificazione}
The Stein fillable contact structure on $M_0$ described by the presentation of
$M$ as $0$-surgery on the right-handed trefoil knot $K$ is $\xi_1$.
\end{prop}
\begin{proof}
By Theorem \ref{honda}, the Stein fillable contact structure on $M_o$ is 
isotopic to $\xi_k$ for some $k \in \N$. 

It is easy to make the meridian of $K$ Legendrian with Thurston--Bennequin 
invariant $-1$ in the standard tight contact structure of $S^3$,
therefore  $tn(F, \xi_k) \geq -1$ because the image of the meridian of $K$ is 
isotopic to $F$ as a framed knot in $M_0$. By Lemma \ref{tn} and Lemma 
\ref{max} this is possible only if $k=1$.
\end{proof}

\subsection{Tight contact structures on $- \Sigma (2,3,6n+5)$} 
The manifold $- \Sigma(2,3,6n+5)$ is  obtained from $M_0$ by $-(n+1)$--surgery on $F$.
For any $n \in \N$ and $n \geq 2$ we define ${\mathcal P}^*_n = \{ -n+1, -n+3, \ldots ,
n-3, n-1 \}$. If $n$ is even, then $0 \notin {\mathcal P}^*_n$ and we define 
${\mathcal P}_n = {\mathcal P}^*_n \cup \{ 0 \}$. In the following we will always 
consider $n$ even, although some of the facts that we are going to prove are 
true for any $n$.

Let $S_+$ and $S_-$ denote the operations of positive and negative stabilisation
defined, for example, in \cite{etnyre:0}, Section 2.7. Given 
$i \in {\mathcal P}^*_n$,
denote the contact structure on $-\Sigma(2,3,6n+5)$ obtained by Legendrian surgery 
on $(M_0, \xi_1)$ along the Legendrian knot $S_+^{(n-1+i)/2} S_-^{(n-1-i)/2}(F)$ by $\eta_i$. 
We denote the tight contact structure on $-\Sigma(2,3,6n+5)$ obtained by Legendrian 
surgery on $(M_0, \xi_n)$ along $F$ by $\eta_0$.

The contact manifolds $(- \Sigma(2,3,6n+5), \eta_i)$ for $i \in {\mathcal P}^*_n$ are the
 Stein fillable contact manifolds considered in \cite{lisca-matic:1}, in fact
$(M_0, \xi_1)$ is the Stein fillable contact manifold obtained by Legendrian 
surgery on a positive trefoil knot in $S^3$ with Thurston-Bennequin invariant 
$0$ by Proposition \ref{identificazione}, and performing Legendrian surgery on
a stabilisation of $F$ is equivalent to performing Legendrian surgery on a 
stabilisation of a meridian of the trefoil knot.
\begin{prop} \label{trucco}
Let $\overline{\eta}_i$ be the contact structure obtained from $\eta_i$ by reversing 
the orientation of the contact planes. Then $\overline{\eta}_i$ is isotopic to
$\eta_{-i}$.
\end{prop}
\begin{proof}
For any $n \in \N^+$ $(M_0, \xi_n)$ is isotopic to $(M_0, \overline{\xi}_n)$. The isotopy 
is induced by a translation in the $t$ direction in the cover $T^2 \times \R$, 
therefore it fixes $F$.
We denote $S_+^{(n-1+i)/2} S_-^{(n-1-i)}(F)$ thought of as a Legendrian knot in
$(M_0, \overline{\xi}_n)$ by $\overline{S_+^{(n-1+i)/2} S_-^{(n-1-i)/2}(F)}$.  Since
changing the orientation of the planes changes positive 
stabilisations into negative ones and vice versa, 
$\overline{S_+^{(n-1+i)/2} S_-^{(n-1-i)/2}(F)}$ is Legendrian isotopic to
$S_+^{(n-1-i)/2} S_-^{(n-1+i)/2}(F)$, therefore inverting the orientation of the planes 
transforms Legendrian surgery on \linebreak $S_+^{(n-1+i)/2} S_-^{(n-1-i)}(F)$ into 
Legendrian surgery on $S_+^{(n-1-i)/2} S_-^{(n-1+i)/2}(F)$.
\end{proof}

\begin{theorem}
The contact structures $\eta_i$ on $- \Sigma(2,3,6n+5)$, with $i \in {\mathcal P}_n$, are 
all pairwise non isotopic. 
\end{theorem}
\begin{proof}
By \cite{lisca-matic:1}, Theorem 4.2, and \cite{lisca-matic:2}, Corollary 4.2,
the contact structures $\eta_i$ with $i \in {\mathcal P}^*_n$ are pairwise non isotopic.
In particular, since we are considering $n$ even, $\eta_i$ is never isotopic to
$\eta_{-i}$ if  $i \in {\mathcal P}^*_n$ because $0 \notin {\mathcal P}^*_n$. Suppose by 
contradiction that $\eta_0$ is 
isotopic to $\eta_i$ for some $i \in {\mathcal P}^*_n$ . Inverting the orientation 
of the contact planes and applying Proposition \ref{trucco}, we obtain
that  $\eta_0$ is also isotopic to $\eta_{-i}$. From this it would follow that $\eta_i$ is 
isotopic to $\eta_{-i}$.
\end{proof}

\begin{rem}
Using methods from \cite{ghiggini-schonenberger} one can prove
that $- \Sigma(2,3,17)$ admits at most three tight contact structures up to 
isotopy, therefore Proposition \ref{trucco} gives the classification of the 
tight contact structures on $- \Sigma (2,3,17)$.
\end{rem}

\subsection{Computation of the homotopy invariants}
In this subsection we will compute the Gompf's three-dimensional homotopy
invariant $d_3(\eta_i)$. This computation will show that all $\eta_i$ are homotopic and
therefore all their Ozsv{\'a}th--Szab{\'o} contact invariants belong to the same 
factor of $\widehat{HF}(-M)$.
 
By \cite{gompf:1}, Theorem 4.5 (for an easy proof of this theorem for integer
homology spheres see also \cite{lisca-matic:1}, Proposition 2.2), $\eta_{i_1}$ is 
homotopic to $\eta_{i_2}$ as a plane field if and only if $d_3(\eta_{i_1})=d_3(\eta_{i_2})$, 
where 
$$d_3 (\eta_i)= \frac 14 (c_1^2(J_i) - 2 \chi(X_i)- 3 \sigma(X_i))$$
and $(X_i, J_i)$ is an almost complex manifold such that $\partial X_i = M$ and 
$\eta_i=TM \cap J(TM)$. 

As almost complex manifold for the computation of $d_3(\eta_i)$ we will
take symplectic fillings of $(M, \eta_i)$ 
endowed with an adapted almost complex structure. More precisely, let $(X_0, \omega)$
 be the weak symplectic filling of $(M_0, \xi_n)$ for any $n \in \N$
constructed in \cite{ding-geiges:1} Proposition 15. If $T \subset M_0$ is a fibre of 
the torus bundle $M_0 \to S^1$, then we can assume that $\int_T \omega =1$. In the setting 
of symplectic fillings Legendrian surgery corresponds to adding symplectic
$2$--handles, so adding symplectic $2$-handles to $(X_0, \omega)$ as explained 
in the definition of $(M, \eta_i)$, we obtain symplectic manifolds
$(X, \omega_i)$ which fill $(M, \eta_i)$ for $i \in {\mathcal P}_n$. We choose almost complex
structures $J_i$ adapted to $\omega_i$ so that the contact structure $\eta_i$ is 
$J_i$--invariant for any $i \in {\mathcal P}_n$, all $J_i$ coincide on $X_0$ and 
the fibre $T$ in $M_0= \partial X_0$ are quasi-complex submanifolds.

In $M_0$, the homology class represented by $F$ is Poincar{\'e} dual of 
$[\omega_0|_{M_0}]$, because $F \cdot T=1= \int_T \omega_0$ and $[T]$ generates $H^2(M_0)$, therefore 
$F$ bounds a surface $\Sigma \subset X_0$ which 
represents the Poincar{\'e} dual of $[\omega]$.
Applying the homology long exact sequence to the pair $(X, X_0)$ we obtain
$H_2(X)=H_2(X_0) \oplus \Z[\overline{\Sigma}]$, where $\overline{\Sigma} \subset X$ is the surface 
obtained by capping $\Sigma$ with the core of the $2$-handle attached along $F^3$.
Analogously, the cohomology exact sequence yields $H^2(X) \cong H^2(X_0) \oplus \Z$, where
the isomorphism is given by $\alpha \mapsto (\iota^* \alpha, \langle \alpha, [\overline{\Sigma}] \rangle)$.

\begin{lemma}\label{serva}
Let $\alpha \in H^2(X)$ be the $2$-dimensional cohomology class determined by $\iota^*(\alpha)=0$ 
and $\langle \alpha, [\overline{\Sigma}] \rangle=1$. Then, up to 
torsion, $\alpha$
is the Poincar{\'e} dual of $[T] \in H_2(X) \cong H_2(X, \partial X)$.
\end{lemma}
\begin{proof}
Any $2$-dimensional homology class can be represented as a closed, oriented 
embedded surface. Let $K$ be a surface representing a homology class in 
$H_2(X_0)$, then $K \cdot T=0$ because $K$ can be made disjoint from $\partial X_0 =M_0$
 and $\langle \alpha, [K] \rangle = \langle \iota^*\alpha, [K] \rangle =0$. On the other hand, $\overline{\Sigma} \cdot T= F \cdot T=
1= \langle \alpha, [\overline{\Sigma}] \rangle$.
\end{proof}
 
\begin{theorem} \label{omotopia}
The contact structures $\eta_i$ with $i \in {\mathcal P}_n$ are pairwise homotopic 
and $d_3(\eta_i)=- \frac 32$.
\end{theorem}
\begin{proof}
To prove that the contact structures are homotopic we will show that 
they have 
the same three dimensional invariant $d_3$. Since in the computation of $d_3(\eta_i)$ 
we use the almost complex manifolds $(X, J_i)$ which are smoothly 
diffeomorphic, it is enough to prove that $c_1^2(J_i)$ does not depend on $i$. 
Given $i_1, i_2 \in {\mathcal P}_n$ we can decompose 
$$c_1^2(J_{i_1})- c_1^2(J_{i_2}) = \langle (c_1(J_{i_1})+ c_1(J_{i_2})), 
PD(c_1(J_{i_1})- c_1(J_{i_2}))\rangle.$$
By the functoriality of the Chern classes for any $i \in {\mathcal P}_n$ we have
$\iota^*(c_1(J_i))=c_1(J_i|_{X_0})$, then $\iota^*(c_1(J_{i_1})- c_1(J_{i_2}))=0$, because all $J_i$ 
 agree on $X_0$. Lemma \ref{serva} implies that $PD(c_1(J_{i_1})- 
c_1(J_{i_2}))$ is a multiple of $[T]$. Since $T$ is a complex submanifold of 
$(X, J_i)$, the adjunction equality gives
$\langle c_1(J_i), [T] \rangle = \chi(T)+ T \cdot T=0$, then $c_1^2(J_{i_1})- c_1^2(J_{i_2}) =0$.

$d_3(\eta_i)$ can be computed for any of the Stein fillable contact structures $\eta_i$ 
with $i \in {\mathcal P}^*_n$ using the Stein filling $(W, J_i)$ described in 
\cite{lisca-matic:1}, Figure 2. One can immediately check that $c_1^2(J_i)=0$, 
$\chi(W)=3$ and $\sigma(W)=0$.
\end{proof}
We stress the point that the Stein manifolds $(W, J_i)$ used to compute $d_3(\eta_i)$
 are different from the almost complex manifolds $(X, J_i)$ used in the first
part of Theorem \ref{omotopia} to show that all $\eta_i$ are homotopic.
\subsection{Computation of the Ozsv{\'a}th-Szab{\'o} invariants}
 In \cite{O-Sz:4}, Section 8, $HF^+( \Sigma(2,3,6n+5))$ is computed. Applying the
long exact sequence relating $HF^+$ and $\widehat{HF}$ and the isomorphism
between $\widehat{HF}_d(Y)$ and $\widehat{HF}_{-d}(-Y)$ it is easy to show that 
$\widehat{HF}(- \Sigma(2,3,6n+5))= (\Z /2 \Z)^{n+1}_{(+2)} \oplus (\Z /2 \Z)^2_{(+1)}$. The degree 
of $c(\xi)$ is $+1$ because $d_3(\eta_i)=- \frac{3}{2}$. By \cite{plam:1}, Section 4 
$\widehat{HF}_{(+1)}(- \Sigma(2,3,6n+5))$ is freely generated by the elements $c(\eta_i)$ 
for $i \in {\mathcal P}^*_n$.

%
%
\noindent
{\em Proof of Theorem \ref{principale}.}
The fix space $Fix(\mathfrak{J}) \subset \widehat{HF}_{(+1)}(- \Sigma(2,3,6n+5))$ is generated
by elements of the form $c(\eta_i)+c(\eta_{-i})$ for $i \in {\mathcal P}^*_n$. 
Let $W$ be the smooth cobordism between $M_0$ and $- \Sigma(2,3,6n+5)$ constructed by 
attaching a $2$-handle to $M_0$ along $F$, then by \cite{O-Sz:cont}, Theorem 4.2
$$\widehat{F}_{\overline W}(c(\eta_i)+c(\eta_{-i}))=\widehat{F}_{\overline W}(c(\eta_i))+ 
\widehat{F}_{\overline W}(c(\eta_{-i}))=2c(\xi_1)=0.$$ 
Consequently $Fix(\mathfrak{J}) \subset \ker \widehat{F}_{\overline W}$, in particular
$$c(\xi_n)= \widehat{F}_{\overline W}(\eta_0) =0$$ because $c(\eta_0) \in Fix(\mathfrak{J})$ by 
Proposition \ref{trucco} and Theorem \ref{coniugazione}. \\
 \qed

In view of Theorem \ref{bruciato} we have the following corollary.
\begin{cor}
The contact manifolds $(M_0, \xi_n)$ are not strongly symplectically fillable
if $n$ is even.
\end{cor}
This is a new non fillability result, because the integer $n_0$ in Theorem 
\ref{geiges} is not given explicitly.

\section{A remark on integer coefficients}
Unfortunately Theorem \ref{principale} does not imply that the 
Ozsv\'ath--Szab\'o contact invariants $c(\xi_n)$ for $n$ even with untwisted integer 
coefficients are zero, but only that they
is the double of some elements of $\widehat{HF}(-M_0) /\pm 1$. Fix an open book 
decomposition of $M_0$ adapted to $\xi_n$ for an even $n$. We denote by $M_0'$ the 
$3$--manifold
obtained by $0$--surgery on the binding and by $M_0''$ the $3$--manifold
obtained by $1$--surgery on the binding. Of course the manifolds $M_0'$ and
$M_0''$ depend on $n$. By \cite{O-Sz:2}, Theorem 9.1 
there is a surgery exact triangle
\begin{align*}
\xymatrix{
\widehat{HF}(-M_0') \ar[rr]^{\widehat{F}} &  & \widehat{HF}(-M_0) \ar[dl] \\
   & \widehat{HF}(-M_0'') \ar[lu] & \\
}
\end{align*}
The group $\widehat{HF}(-M_0')$ is generated by $c_0$, therefore if 
$F(c_0)= c(\xi_2) \neq 0$, the exact triangle becomes a short exact sequence
$$
\begin{matrix}
0 & \to & \widehat{HF}(-M_0') & \to & \widehat{HF}(-M_0) & \to & \widehat{HF}(-M_0'') 
& \to & 0
\end{matrix}$$
If $c(\xi_n)$ is non primitive there are torsion elements in 
$\widehat{HF}(-M_0'')$. Since all Heegaard--Floer homology groups known so far
are free, it is reasonable to expect that $c(\xi_n)=0$ also in the 
Heegaard--Floer homology group with integer coefficients. 

\bibliographystyle{amsplain}
\bibliography{contatto}
\end{document}